\documentclass[12pt]{article}
\usepackage{latexsym}
\usepackage{multirow}
\usepackage{graphicx}
\usepackage{float}
\usepackage{authblk}%%%for authors
\usepackage{multirow}
\usepackage{kantlipsum}
\usepackage{setspace}
\usepackage{xcolor, color}
\usepackage{bbm}  %示性函数
\def\q{\hfill\rule{1ex}{1ex}}
\def\0{\emptyset}

\def\p{{\bf Proof.~~}}

\usepackage{indentfirst}
\usepackage{amsmath}

\def\0{\emptyset}

\def\p{{\bf Proof.} \quad}
\def\q{\hfill\rule{1ex}{1ex}}

\usepackage{mathrsfs}
\usepackage{amssymb}
\usepackage{listings}
\usepackage{amsmath}
\usepackage{pict2e}
\usepackage{amsfonts}
\usepackage{graphicx}
\usepackage[small]{caption}
\usepackage{subfigure}
\usepackage{cite}
\usepackage{mathtools}
\usepackage{multirow}
\usepackage{tikz}
\newtheorem{theorem}{Theorem}[section]
\newtheorem{corollary}[theorem]{Corollary}

	\newtheorem{Lemma}[theorem]{Lemma}

\newtheorem{proposition}[theorem]{Proposition}

\usepackage{graphicx}
\newcommand{\medcup}{\mathbin{\scalebox{1.3}{$\cup$}}}
\newcounter{countclaim}

\begin{document}
	\title{Spectral radius and homeomorphically
		irreducible spanning trees of graphs
		\thanks{This research was partially supported by NSFC (No. 12371345 and 12471325). }}
	
	\author[1,2]{\small Bingqian Gao\thanks{Email: bqgao\_0110@163.com}}
	
	\author[1,2]{\small Huiqing Liu\thanks{Email: hqliu@hubu.edu.cn}}
	
	\author[3]{\small Jing Zhao\thanks{Email: jingzhao\_0911@163.com}}
	
	\affil[1]{\footnotesize Hubei Key Laboratory of Applied Mathematics, Faculty of Mathematics and Statistics, 
		
		Hubei University, Wuhan 430062, China}
	
	\affil[2]{\footnotesize Key Laboratory of Intelligent Sensing System and Security (Hubei University), Ministry of Education}
	
	\affil[3]{\footnotesize School of Mathematics and Physics, Anqing Normal University, Anqing 246133, China}

	\date{}
	
	\maketitle
	
	\begin{abstract}
		\baselineskip=0.5cm
		For a connected graph $G$, a spanning tree $T$ of $G$ is called a homeomorphically irreducible spanning tree (HIST) if $T$ has no vertices of degree 2. Albertson {\em et al.} proved that it is $NP$-complete to decide whether a graph contains a HIST. In this paper, we provide some spectral conditions that guarantee the existence of a HIST in  a connected graph. Furthermore, we also present some sufficient conditions in terms of the order of a graph $G$ to ensure the existence of a HIST in $G$.

		\vskip 0.2cm
		
		{\bf AMS} classification: 05C35, 05C50
		
		\vskip.2cm
		
	\end{abstract}
	\vskip 0.2cm

	\section{Introduction}
	Throughout this paper, we consider only simple, undirected and connected graphs. Let $G=(V(G), E(G))$ be a graph. For $v\in V(G)$, denoted by $N_G(v)$ the set of all vertices in $G$ adjacent to $v$, and $N_G[v]=N_G(v)\cup\{v\}$. Set $d_G(v)=|N_G(v)|$, namely, the {\em degree} of the vertex $v$ in $G$, where the index $G$ will be omitted if there is no risk of confusion. Let $\Delta(G)$ and $\delta(G)$ be the {\em maximum degree} and {\em minimum degree} of $G$, respectively. For $S\subseteq V(G)$, we denote by $G[S]$ the subgraph of $G$ induced by $S$, and $G-S:=G[V(G)\setminus S]$. If $S=\{v\}$, then we simplify $G-\{v\}$ to $G-v$. Given two disjoint vertex sets $X$ and $Y$ of $G$, let $E(X,Y)$ be the set of all edges with one end in $X$ and the other end in $Y$. If $X=\{x\}$, we simply write $E(x,Y)$ for $E(X,Y)$. A nontrivial graph $G$ is said to be {\em $k$-connected} if there does not exist a vertex cut of size $k-1$ whose removal disconnects the graph. We will use $P_n$, $K_n$ to denote a path, a complete graph of order $n$, respectively, and use $K_{p,q}$ to denote a complete bipartite graph with bipartition $(S,T)$, where $|S|=p$ and $|T|=q$.
	
	Let $A(G)=(a_{ij})$ be the \emph{adjacency matrix} of $G$, where $a_{ij}=1$ if $v_i$ is adjacent to $v_j$, and $a_{ij}=0$ otherwise, where $v_i,v_j\in V(G)$. It follows immediately that if $G$ is a simple graph, then $A(G)$ is a symmetric $(0,1)$ matrix in which every diagonal entry is zero. Since $A(G)$ is real and symmetric, its eigenvalues are real. The {\em spectral radius} of $G$, denoted by $\rho(G)$, is the largest eigenvalue of $A(G)$. Note that if $G$ is connected, then $A(G)$ is irreducible, and so by the Perron-Frobenius theory of non-negative matrices, $\rho(G)$ has multiplicity one and there exists a unique positive unit eigenvector (also called {\em Perron-eigenvector}) corresponding to $\rho(G)$.
	
	Given a graph $F$, we say that $G$ is {\em $F$-free} if it does not contain $F$ as a subgraph. In 2010, Nikiforov \cite{Niki}  proposed a spectral extremal problem as follows.
	
	\vskip.2cm
	
	\noindent{\bf Problem 1.} For a given graph $F$, what is the maximum spectral radius of an $n$-vertex $F$-free graph?
	
	\vskip.2cm
	
	Recently, much attention has been paid to the various families of graphs $F$, such as $K_r$ \cite{Niki2007,Wilf}, $K_{s,t}$ \cite{Babai2009}, $K_{s,t}$-minor \cite{Zhai2022J}, cycles \cite{Lu2024,Min2022,Niki,Li2022,Zhai2012,Zhai2020}, friendship graph \cite{Zhai2022E}, linear forest \cite{Chen2019} and references therein. 
	
	For a connected graph $G$, a spanning tree $T$ of $G$ is called a {\em homeomorphically irreducible spanning tree} (HIST) if $T$ has no vertices of degree 2. Intuitively, HIST represents a class of graphs that are, in a sense, antithetical to Hamiltonian paths. In \cite{Michael1990}, Albertson {\em et al.} showed that it is $NP$-complete to decide if a graph has a HIST. For more results related to HIST, please refer to \cite{Furuya2025,Hill1974,Ito2022,Li2024}. In this paper, we aim to provide some spectral conditions to guarantee the existence of a HIST in connected graphs.

	We are now going to introduce two graphs $L_n$ and $B_n$ of order $n$ shown in Figure 1, before presenting our main results in this paper. Let $L_n$ be a graph obtained from $K_{n-2}$ and $P_2$ by joining one pendant vertex of $P_2$ to one vertex of $K_{n-2}$ (see Figure 1 (a)). Let $B_n$ be a graph obtained from $K_{n-3}$ and $P_3$ by joining two pendant vertices of $P_3$ to two distinct vertices of $K_{n-3}$, respectively (see Figure 1 (b)).
	
	\begin{center} \setlength{\unitlength}{1.2mm}
		\begin{picture}(150,40)
			\thicklines % 使线条更粗
			\put(25, 30){\circle{20}}
			
			\put(35, 30){\line(1, 0){10}}
			\put(45, 30){\line(1, 0){10}} 
			
			\put(35, 30){\circle*{1.5}}
			\put(45, 30){\circle*{1.5}}
			\put(55, 30){\circle*{1.5}}

			%\put(46, 31){$v$}
			
			\put(22, 29){$K_{n-2}$}
			
			%\put(36.5,32){\makebox(0,0){$v_3$}}
			%\put(46,32){\makebox(0,0){$v_2$}}
			%\put(56,32){\makebox(0,0){$v_1$}}
			
			%%%%%%%%%%%%%%%%%%%%%%%%%%%%%%%%%%%%%%%%%%%%%%%%%%		
			
			\put(90, 30){\circle{20}} 
			
			\put(90, 40){\circle*{1.5}}
			\put(110, 40){\circle*{1.5}}
			\put(110, 30){\circle*{1.5}}
			\put(90, 20){\circle*{1.5}}
			\put(110, 20){\circle*{1.5}}
			
			\put(90, 40){\line(1, 0){20}}
			\put(90, 20){\line(1, 0){20}} 
			
			\put(110, 40){\line(0, -1){20}}
			
			%\put(90, 42){$u_1$}\put(90, 17){$u_2$}\put(110, 42){$v_1$}\put(110, 17){$v_2$}

			\put(87.5, 29){$K_{n-3}$}
			
			%\put(90,42){\makebox(0,0){$u_5$}}
			%\put(90,17){\makebox(0,0){$u_4$}}
			%\put(113,40){\makebox(0,0){$u_1$}}
			%\put(113,30){\makebox(0,0){$u_2$}}
			%\put(113,20){\makebox(0,0){$u_3$}}
			
			%%%%%%%%%%%%%%%%%%%%%%%%%%%%%%%%%%%%%%%%%%%%%%%%%
			
			\put(37.5, 10){$(a)$}\put(88.5, 10){$(b)$}
			
			\put(25, 0){Figure 1. (a) The graph $L_n$, and (b) the graph $B_n$}
		\end{picture}
	\end{center}

	\begin{theorem}\label{Thm1}
		Let $G$ be a connected graph of order $n\ge 7$. If $\rho(G)\ge\rho(L_n)$, then $G$ contains a HIST unless $G \cong L_n$.
	\end{theorem}
	
	We also provide a spectral condition to guarantee the existence of a HIST in a $2$-connected graph.
	
	\begin{theorem} \label{Thm2}
		Let $G$ be a 2-connected graph of order $n\ge 8$. If $\rho(G)\ge\rho(B_n)$, then $G$ contains a HIST unless $G \cong B_n$.
	\end{theorem}

	\noindent{\bf Remark 1. }	The graphs $L_n$ and $B_n$ does not contain HISTs. The details are given in Proposition \ref{le00}.\medskip
	
	The above two theorems can be regarded as answers to Problem $1$ when $G$ is HIST-free:
	
	\begin{itemize}
		\item For every connected HIST-free graph $G$ of order $n\ge 7$, we have $\rho(G)\le\rho(L_n)$ with equality if and only if $G \cong L_n$.
		\item For every $2$-connected HIST-free graph $G$ of order $n\ge 8$, we have $\rho(G)\le\rho(B_n)$ with equality if and only if $G \cong B_n$.  
	\end{itemize}
	
	By Theorems \ref{Thm1} and \ref{Thm2}, we can obtain the following sufficient conditions in terms of the order of $G$ for the existence of a HIST in $G$. %The details are given in Section 3.
	
	\begin{corollary}\label{Thm3}
		Let $G$ be a connected graph of order $n\ge 7$. If $\rho(G)\ge n-3+\frac{1}{n-3}$, then $G$ contains a HIST.
	\end{corollary}
	
	\begin{corollary} \label{Thm4}
		Let $G$ be a 2-connected graph of order $n\ge 8$. If $\rho(G)\ge n-4+\frac{2}{n-4}$, then $G$ contains a HIST.
	\end{corollary}

	%-------------	
	\section{Preliminaries}
	
	In this section, we first explain why the conclusion of Theorem \ref{Thm1} excludes the graph $L_n$ and the conclusion of Theorem \ref{Thm2} excludes the graph $B_n$.	
	%------------	
	\begin{proposition}\cite{Li2024} \label{le0}
		Any connected graph with a 
		cut-vertex of degree $2$ 
		has no HISTs. 
	\end{proposition}	
	
	%--------------	
	\begin{proposition} \label{le00}
		Let $G$ be a connected graph. If $G$ contains a path $P_5=s_0s_1s_2s_3s_4$ with $d(s_0),d(s_4)\ge3$, $d(s_1)=d(s_2)=d(s_3)=2$, then $G$ does not contain HISTs.
	\end{proposition}
	\p	
	Suppose, to the contrary, that $G$ has a HIST $T$. Then $\{s_1s_2, s_2s_3\}\cap E(T)\not=\emptyset$ as $d(s_2)=2$. Without loss of generality, we assume $s_1s_2\in E(T)$. Since $T$ does not contain vertex of degree 2, we have $s_2s_3\notin E(T)$. Similarly, $s_0s_1\notin E(T)$. Hence, $G[\{s_1,s_2\}]$ is a component of $T$, a contradiction with the fact that $T$ is connected.
	\q \medskip
	
	The following conclusion then follows directly from Propositions \ref{le0} and \ref{le00}.
	
	\begin{corollary}\label{cor1}
		The graphs $L_n$ and $B_n$ do not contain HISTs. 
	\end{corollary}
	
	We use the following lemma to deal with the spectral radii of subgraphs. 
	
	\begin{Lemma}\cite{CRS10} \label{le1}
		Let $G'$ be a proper subgraph of $G$. Then $\rho(G')<\rho(G)$. 
	\end{Lemma}
	
	Next, we will present some upper bounds of the spectral radii of $L_n$ and $B_n$.
	
	\begin{proposition} \label{2.6}
		(i) The spectral radius of the graph $L_n$ is the largest root of the equation 
		$x^4-(n-4)x^3-(n-1)x^2+(2n-8)x+n-3=0$. 
		
		(ii)  The spectral radius of the graph $B_n$ is the largest root of the equation 
		$x^4-(n-5)x^3-(n-1)x^2+(3n-16)x+2n-8=0$.
		
		Moreover, $\rho(L_n)<n-3+\frac{1}{n-3}$ and $\rho(B_n)< n-4+\frac{1}{n-4}$.	
	\end{proposition}
	\p (i)  Let $\lambda$ be the spectral radius of $L_n$. We label the vertices of $L_n$ by $\{v_1,v_2,\ldots,v_n\}$ such that $L_n[\{v_3, v_4,\ldots,v_n\}] \cong K_{n-2}$ and $L_n[\{v_1,v_2,v_3\}] \cong P_3$, where $v_2$ is the vertex of degree 2. Let $\bf{x}$$=(x_1, x_2,\ldots, x_n)^t$ be the Perron-eigenvector of $L_n$ with coordinate $x_i$ corresponding to vertex $v_i$. From $\lambda \bf{x}$$=A(L_n)\bf{x}$, we have $x_4=\cdots=x_n$. Then
	\[
	\begin{cases} 
		\lambda x_1=x_2, \\
		\lambda x_2=x_1+x_3, \\ 
		\lambda x_3=x_2+(n-3)x_4, \\ 
		\lambda x_4=x_3+(n-4)x_4. 
	\end{cases}
	\]
	Hence $\lambda^4-(n-4)\lambda^3-(n-1)\lambda^2+(2n-8)\lambda+n-3=0$.
	
	Since $K_{n-2}$ is a proper subgraph of $L_n$, we obtain $\lambda>n-3$ by Lemma \ref{le1}, and then, we let $\lambda=n-3+t$, where $t>0$. Note that $\lambda$ is the largest root of the equation 
	$$x^4-(n-4)x^3-(n-1)x^2+(2n-8)x+n-3=0.$$ This implies 
	\begin{eqnarray}\label{31}
		t^4+(3n-8)t^3+(3n^2-16n+19)t^2+(n^3-8n^2+19n-14)t+3-n=0.
	\end{eqnarray}
	Note that $t>0$ and
	\begin{eqnarray}\label{32}
		\left\{
		\begin{array}{l}
			3n-8>0, \\ 
			3n^2-16n+19=n(3n-16)+19>0, \\ 
			n^3-8n^2+19n-14>(n-3)^2>0,
		\end{array}
		\right.
	\end{eqnarray}
	as $n\ge 7$. Therefore, by \eqref{31} and \eqref{32}, we have that
	$(n^3-8n^2+19n-14)t+3-n<0,$
	it follows that
	$$t<\frac{n-3}{n^3-8n^2+19n-14}<\frac{1}{n-3},$$
	and then 
	$$\lambda=n-3+t<n-3+\frac{1}{n-3}.$$
	
	%By Theorem \ref{Thm1}, we have that if $\rho(G)\ge n-3+\frac{1}{n-3}$, then $G$ contains a HIST. Therefore, the proof of Corollary \ref{Thm3} is complete.\q	

	(ii) Let $\rho$ be  the spectral radius of $B_n$.  We label the vertices of $B_n$ by $\{u_1,u_2,\ldots,u_n\}$ such that $B_n[\{u_1,\ldots,u_5\}] \cong C_5$ and $B_n[\{u_4,u_5,\ldots,u_n\}] \cong K_{n-3}$. Let $\bf{y}$$=(y_1, y_2,\ldots, y_n)^t$ be the Perron-eigenvector of $A(B_n)$ with coordinate $y_i$ corresponding to vertex $u_i$. From $\rho \bf{y}$$=A(B_n)\bf{y}$, we have $y_1=y_3$, $y_4=y_5$ and $y_6=\cdots=y_n$. Then
	\[
	\begin{cases} 	
		\rho y_1=y_2+y_5, \\ 
		\rho y_2=2y_1, \\ 
		\rho y_5=y_1+y_5+(n-5)y_6, \\ 
		\rho y_6=2y_5+(n-6)y_6. 
	\end{cases}
	\]
	Hence $\rho^4-(n-5)\rho^3-(n-1)\rho^2+(3n-16)\rho+2n-8=0$.\vskip.2cm
	
	Note that $\rho\ge \rho(K_{n-3})>n-4$ by Lemma \ref{le1}, and then, we can suppose $\rho=n-4+s$, where $s>0$.  Hence
	\begin{eqnarray}\label{41}
		s^4+(3n-11)s^3+(3n^2-22n+37)s^2+(n^3-11n^2+37n-40)s+8-2n=0.
	\end{eqnarray}
	Note that $s>0$ and 
	\begin{eqnarray}\label{42}
		\left\{
		\begin{array}{l}
			3n-11>0, \\ 
			3n^2-22n+37=n(3n-22)+37>0, \\ 
			n^3-11n^2+37n-40>(n-4)^2>0,
		\end{array}
		\right.
	\end{eqnarray}
	as $n\ge 8$. Hence, by \eqref{41} and \eqref{42}, we obtain that
	$(n^3-11n^2+37n-40)s+8-2n<0,$
	it follows that
	$$s<\frac{2n-8}{n^3-11n^2+37n-40}<\frac{2}{n-4},$$
	then 
	$$\rho=n-4+s<n-4+\frac{2}{n-4}.$$
	
	Therefore the proof of Proposition \ref{2.6} is complete.	\q	\medskip
	
	\noindent{\bf Remark 2. }	Corollaries \ref{Thm3} and \ref{Thm4} follow directly from  Proposition \ref{2.6}, Theorems \ref{Thm1} and \ref{Thm2}, respectively.	\medskip
	%---
	
	In the remainder of this section, we present some results that will be used in the following sections.

	\begin{Lemma} \cite{CRS10} \label{le2} 
		Let $G$ be a graph with maximum degree $\Delta(G)$. Then $\rho(G)\leq \Delta(G)$.
	\end{Lemma}
	
	\begin{Lemma} \cite{Hong1988,Hong1999} \label{le3}
		Let $G$ be a simple connected graph of order $n$ and size $m$, and let $\delta=\delta(G)$. Then
		$$\rho (G)\le \frac{\delta-1+\sqrt{(\delta+1)^2+4(2m-\delta n)}}{2}.$$ 
		Moreover, if $\delta =1$, then $\rho(G)\leq\sqrt{2m-n+1}$.
		
	\end{Lemma}
	
	%------------

	%%%%%%%%%%%%	
	\section{Proofs}
	In this section, we will prove Theorems \ref{Thm1} and \ref{Thm2}.
	
	\vskip.2cm	
	
	\noindent{\bf Proof of Theorem \ref{Thm1}.}
	First, we note that $K_{n-2}$ is a proper subgraph of $L_n$, then 
	\begin{eqnarray} \label{q1}
		\rho(L_n)>\rho(K_{n-2})=n-3.
	\end{eqnarray}
	By Lemma \ref{le2}, we have $\rho(G)\le\Delta(G)$. Hence $\Delta(G)\ge\rho(G)\ge\rho(L_n)>n-3$, that is, $\Delta(G)\ge{n-2}.$
	
	Let $x\in{V(G)}$ with $d(x)=\Delta(G)$. Denote $N(x)=\{x_1,x_2,\ldots,x_{\Delta(G)}\}$. If $\Delta(G)=n-1$, then $G$ has a HIST $T$ with $E(T)=E(x,N(x))$. So, in the following, we assume that $\Delta(G)=n-2$.
	
	Denote $V(G)\backslash N[x]=\left\{y\right\}$. Because $G$ is connected, $E(y,N(x))\neq \emptyset$. Without loss of generality, we assume $N(y)\cap N(x)=\{x_1,x_2,\ldots,x_a\}$, where $1\leq a \le n-2$. If there exists some $x_i ~(1\leq i\leq a)$ such that $N(x_i)\cap N(x)\neq \emptyset$, say $x_ix_j$ for some $j\neq i$, then $G$ has a HIST $T$ with $E(T)=E(x,N(x)\backslash \{x_j\})\bigcup\{x_iy,x_ix_j\}$. So $N(x_i)\cap N(x)=\emptyset$ for all $1\le i\le a$, that is $d(x_i)=2$ for $1\le i\le a$. \medskip
	
	\noindent{\bf Fact 1.} $a=1$.\medskip
	
	\noindent{\bf Proof of Fact 1.}	Suppose that $a\geq2$. If $a=n-2$, then $|E(G)|=2(n-2)$, and then by Lemma \ref{le3} and \eqref{q1}, $\rho(G)\leq \sqrt{2|E(G)|-n+1} =\sqrt{3n-7} <\sqrt{n^2-6n+9}= n-3<\rho(L_n)$ as $n\ge 7$, a contradiction with the assumption. Therefore $2\le a\le n-3$. Note that $N(x_j)\subseteq \{x, x_{a+1},\ldots, x_{n-2}\}$, and hence $d(x_j)\le n-a-2$ for all $a+1\le j\le n-2$. Thus
	\begin{eqnarray*}
		2|E(G)|&=&d(x)+d(y)+\sum_{i=1}^{a}d(x_i)+\sum_{j=a+1}^{n-2}d(x_j)\\
		&\le &(n-2)+a+2a+(n-a-2)^{2}=a^2+(7-2n)a+n^2-3n+2\\
		&\le& \max \{2^2+2(7-2n)+n^2-3n+2,(n-3)^2+(n-3)(7-2n)+n^2-3n+2\}\\
		&=&n^2-7n+20.
	\end{eqnarray*}
	By Lemma \ref{le3}, together with \eqref{q1} and $n\ge 7$, one obtains 
	\begin{eqnarray*}
		\rho(G)\leq \sqrt{2|E(G)|-n+1}
		\le\sqrt{n^2-8n+21}
		\le\sqrt{n^2-6n+9}
		= n-3
		<\rho(L_n),
	\end{eqnarray*}
	a contradiction with the assumption. \q
	
	By Fact 1, $a=1$. Recall that $d(x_1)=2$. Therefore $G$ is a subgraph of $L_n$. Furthermore, $G\cong L_n$, for otherwise, if $G$ is a proper subgraph of $L_n$, then by Lemma \ref{le1}, $\rho(G)<\rho(L_n)$, a contradiction.
	
	Therefore, the proof of Theorem \ref{Thm1} is complete. \q
	
	\vskip.2cm
	
	Next, we are going to give a proof of Theorem \ref{Thm2}.	
	
	\vskip.2cm
	
	\noindent{\bf Proof of Theorem \ref{Thm2}.} First, we note that 
	\begin{eqnarray}\label{q2}
		\rho(B_n)>\rho(K_{n-3})=n-4.
	\end{eqnarray} 

	Let $G$ be a $2$-connected graph of order $n\ge 8$ satisfying $\rho(G)\ge \rho(B_n)$. Then $\delta(G)\ge 2$, and then by Lemma \ref{le3}  and \eqref{q2}, we have $$n-4<\rho(B_n)\le \rho(G)\le \frac{1+\sqrt{9+4(2|E(G)|-2 n)}}{2},$$ that is, 
	\begin{eqnarray} \label{q4}
		2|E(G)|>\frac{(2n-9)^2+8n-9}{4}=n^2-7n+18.
	\end{eqnarray}
	
	By Lemma \ref{le2}, we have $\Delta(G)\ge\rho(G)\ge\rho(B_n)>n-4$, that is, $\Delta(G)\ge n-3$.  Let $u\in V(G)$ with $d(u)=\Delta(G)$. Denote $N(u)=\{u_1,u_2,\ldots,u_{\Delta(G)}\}$.
	If $\Delta(G)=n-1$, then $G$ contains a HIST $T$ with $E(T)=E(u,N(u))$. So, in the following, we consider two cases.
	
	\vskip.2cm
	
	\noindent {\bf Case 1.}
	$\Delta(G)=n-2$.
	
	\vskip.2cm
	
	In this case, we let $V(G)\backslash N[u]=\{v\}$. Denote $N(v)=\{u_1,u_2,\ldots,u_b\}$. Since $G$ is 2-connected, $2 \le b \le n-2$. First we can assume that $N(v)$ is an independent set, for otherwise, there exist some $u_r$ and $u_s$ $(1\le r,s \le b)$ such that $u_ru_s\in E(G)$, then $G$ has a HIST $T$ with $E(T)=E(u,N(u)\backslash \{u_s\})\bigcup \{u_rv,u_ru_s\}$. \medskip
	
	\noindent{\bf Fact 2.} $b< n-2$.\medskip
	
	\noindent{\bf Proof of Fact 2.}	Suppose that $b= n-2$. Then $N(u)=N(v)$ is an independent set, furthermore, $G\cong K_{2,n-2}$, and thus, $2|E(G)|=4(n-2)<n^2-7n+18$ as $n\ge 8$, a contradiction with \eqref{q4}. \q \medskip
	
	By Fact 2, $N(u)\setminus N(v)\neq \emptyset$. On the other hand, since $G-u$ is connected, we have $E(N(v), N(u)\setminus N(v))\neq \emptyset$, which implies that there exist some $u_i ~(1\leq i\leq b)$ and $u_j$ (${b+1}\leq j\leq \Delta(G)$) such that $u_iu_j\in{E(G)}$. Then $G$ contains a HIST $T$ with $E(T)=E(u,N(u)\backslash \{u_j\})\bigcup\{vu_i,u_iu_j\}$.
	
	\vskip.3cm
	
	\noindent {\bf Case 2.} 
	$\Delta(G)=n-3$.
	
	\vskip.2cm
	
	In this case, we let $v_1,v_2\in V(G)\backslash N[u]$. If $N(v_1)\cap N(v_2)\neq\emptyset$, say $u_1\in N(v_1)\cap N(v_2)$, then $G$ has a HIST $T$ with $E(T)=E(u,N(u))\bigcup\{u_1v_1,u_1v_2\}$. So, in the following, we can assume $N(v_1)\cap N(v_2)=\emptyset$. We consider two subcases. 
	
	\vskip.2cm
	
	\noindent {\bf Subcase 2.1.} $v_1v_2\in E(G)$.
	
	\vskip.2cm
	
	In this subcase, we let $N(v_1) \backslash \{v_2\}=\{u_1,u_2,\ldots,u_c\}$ and $N(v_2) \backslash \{v_1\}=\{u_{c+1},\ldots,u_d\}$, where $c\ge1$, $d-c\ge1$ and $d\le n-3$. Denote $X_1:=N(v_1) \backslash \{v_2\}$, $X_2:=N(v_2) \backslash \{v_1\}$ and $X:=N(u)\setminus(N(v_1)\cup N(v_2))$.  \medskip
	
	\vskip.2cm
	
	\begin{center}
		\begin{tikzpicture}[line width =0.9pt]			
			
			\filldraw (1.5,-1.8) circle (.07);	
			\filldraw (0,-3) circle (.07);	
			\filldraw (1,-3) circle (.07);	
			
			\filldraw (2,-3) circle (.07);
			\filldraw (2.28,-3) circle (.01);		
			\filldraw (2.5,-3) circle (.01);
			\filldraw (2.72,-3) circle (.01);
			\filldraw (3,-3) circle (.07);	
			
			\filldraw (0.5,-4.2) circle (.07);	
			\filldraw (2.5,-4.2) circle (.07);	
			
			\draw (1.5,-1.8) -- (0,-3);
			\draw [red](1.5,-1.8) -- (1,-3);
			\draw (1.5,-1.8) -- (2,-3);
			\draw [red](1.5,-1.8) -- (3,-3);
			
			\draw(0.5,-4.2) -- (0,-3);
			\draw[red] (2.5,-4.2) -- (1,-3);
			
			\draw [red](2.5,-4.2) -- (2,-3);
			\draw (2.5,-4.2) -- (3,-3);
			
			\draw [red](0.5,-4.2) -- (2.5,-4.2);
			\draw [red](0,-3) to [bend left=20] (1,-3); 
			
			\node at(1.2,-1.8) {$u$};
			\node at(-0.3,-3) {$u_1$};
			\node at(0.7,-3.2) {$u_2$};
			\node at(1.7,-3) {$u_3$};
			\node at(3.3,-3) {$u_d$};
			\node at(0.5,-4.5) {$v_1$};
			\node at(2.5,-4.5) {$v_2$};
			
			\node at(1.5,-5.5) {$(a)~c=1$};
			
			%%%%%%%%%%%%%%%%%%%%%%%%%%%%%%%%%%%%%%%%%%%	
			
			\filldraw (6.5,-1.8) circle (.07);
			
			\filldraw (5,-3) circle (.07);
			\filldraw (5.28,-3) circle (.01);		
			\filldraw (5.5,-3) circle (.01);
			\filldraw (5.72,-3) circle (.01);
			\filldraw (6,-3) circle (.07);
			
			\filldraw (7,-3) circle (.07);
			\filldraw (7.28,-3) circle (.01);		
			\filldraw (7.5,-3) circle (.01);
			\filldraw (7.72,-3) circle (.01);
			\filldraw (8,-3) circle (.07);
			
			\filldraw (5.5,-4.2) circle (.07);
			\filldraw (7.7,-4.2) circle (.07);
			
			\draw (6.5,-1.8) -- (5,-3);
			\draw[red](6.5,-1.8) -- (6,-3);
			\draw (6.5,-1.8) -- (7,-3);
			\draw [red](6.5,-1.8) -- (8,-3);
			
			\draw [red](5.5,-4.2) -- (5,-3);
			\draw [red](5.5,-4.2) -- (6,-3);
			\draw (7.7,-4.2) -- (7,-3);
			\draw (7.7,-4.2) -- (8,-3);

			\draw [red](5.5,-4.2) -- (7.7,-4.2);
			
			\draw[red] (6,-3) to [bend left=20] (7,-3); 
			
			\node at(6.2,-1.8) {$u$};
			\node at(4.7,-3) {$u_1$};
			\node at(5.8,-2.8) {$u_c$};
			\node at(7.3,-2.8) {$u_{c+1}$};
			\node at(8.3,-3) {$u_d$};
			\node at(5.5,-4.5) {$v_1$};
			\node at(7.77,-4.5) {$v_2$};
			
			\node at(6.5,-5.5) {$(b)~c\ge2$};
			
			%%%%%%%%%%%%%%%%%%%%%%%%%%%%%%%%%%%%%%%%%%%		
			
		\end{tikzpicture}
	\end{center}	  
	
	\vskip.2cm
	\begin{center}
		\begin{tikzpicture}[line width =0.9pt]	
			
			\filldraw (4.5,-1.8) circle (.07);
			
			\filldraw (2,-3) circle (.07);
			\filldraw (3,-3) circle (.07);
			\filldraw (3.28,-3) circle (.01);		
			\filldraw (3.5,-3) circle (.01);
			\filldraw (3.72,-3) circle (.01);
			\filldraw (4,-3) circle (.07);
			
			\filldraw (5,-3) circle (.07);
			\filldraw (5.28,-3) circle (.01);		
			\filldraw (5.5,-3) circle (.01);
			\filldraw (5.72,-3) circle (.01);
			\filldraw (6,-3) circle (.07);
			\filldraw (7,-3) circle (.07);
			
			\filldraw (3,-4.2) circle (.07); 
			\filldraw (6,-4.2) circle (.07); 
			
			\draw [red](4.5,-1.8) -- (2,-3);
			\draw (4.5,-1.8) -- (3,-3);
			\draw [red](4.5,-1.8) -- (4,-3);
			\draw [red](4.5,-1.8) -- (5,-3);
			\draw (4.5,-1.8) -- (6,-3);
			\draw [red](4.5,-1.8) -- (7,-3);
			
			\draw[red] (3,-4.2) -- (2,-3);
			\draw (3,-4.2) -- (3,-3);
			\draw (3,-4.2) -- (4,-3);
			
			\draw (6,-4.2) -- (5,-3);
			\draw (6,-4.2) -- (6,-3);
			\draw [red](6,-4.2) -- (7,-3);
			
			\draw (3,-4.2) -- (6,-4.2);
			
			\draw[red] (2,-3) to [bend left=20] (3,-3); 
			\draw[red] (6,-3) to [bend left=20] (7,-3); 
			
			\node at(4.2,-1.7) {$u$};
			\node at(1.8,-2.8) {$u_1$};
			\node at(2.9,-2.8) {$u_2$};
			\node at(3.8,-2.8) {$u_c$};
			\node at(4.7,-3.2) {$u_{c+1}$};
			\node at(6.4,-3.2) {$u_{d-1}$};
			\node at(7.3,-3.2) {$u_d$};
			\node at(3,-4.5) {$v_1$};
			\node at(6,-4.5) {$v_2$};

			%%%%%%%%%%%%%%%%%%%%%%%%%%%%%%%%%%%%%%%%%%
			
			\filldraw (11.5,-1.8) circle (.07);
			
			\filldraw (9,-3) circle (.07); 
			\filldraw (10,-3) circle (.07);
			\filldraw (10.28,-3) circle (.01);		
			\filldraw (10.5,-3) circle (.01);
			\filldraw (10.72,-3) circle (.01);
			\filldraw (11,-3) circle (.07); 
			
			\filldraw (12,-3) circle (.07);  
			\filldraw (12.28,-3) circle (.01);		
			\filldraw (12.5,-3) circle (.01);
			\filldraw (12.72,-3) circle (.01);
			\filldraw (12.9,-3) circle (.07); 
			\filldraw (14,-3) circle (.07); 
			
			\filldraw (10,-4.2) circle (.07); 
			\filldraw (12.5,-4.2) circle (.07); 
			
			\draw [red](11.5,-1.8) -- (9,-3); 
			\draw (11.5,-1.8) -- (10,-3);
			\draw (11.5,-1.8) -- (11,-3);
			\draw [red] (11.5,-1.8) -- (12,-3);
			\draw [red](11.5,-1.8) -- (12.9,-3);
			\draw [red] (11.5,-1.8) -- (14,-3);
			
			\draw [red] (10,-4.2) -- (9,-3);
			\draw (10,-4.2) -- (10,-3);
			\draw [red](10,-4.2) -- (11,-3);
			
			\draw (12.5,-4.2) -- (12,-3);
			\draw (12.5,-4.2) -- (12.9,-3);
			\draw (12.5,-4.2) -- (14,-3);
			
			\draw[red] (10,-4.2) -- (12.5,-4.2);
			
			\draw[red] (9,-3) to [bend left=20] (10,-3); 	
			
			\node at(11.2,-1.7) {$u$};
			\node at(8.9,-2.8) {$u_1$};
			\node at(9.9,-2.8) {$u_2$};
			\node at(10.8,-2.8) {$u_c$};
			\node at(11.7,-3.2) {$u_{c+1}$};
			\node at(13.3,-3.2) {$u_{d-1}$};
			\node at(14.3,-3.2) {$u_d$};
			\node at(10,-4.5) {$v_1$};
			\node at(12.5,-4.5) {$v_2$};
			
			\node at(4.5,-5.5) {$(c)$~neither $X_1$ nor $X_2$ is independent};
			\node at(11.5,-5.5) {$(d)~|X_1|\ge 3$};
			
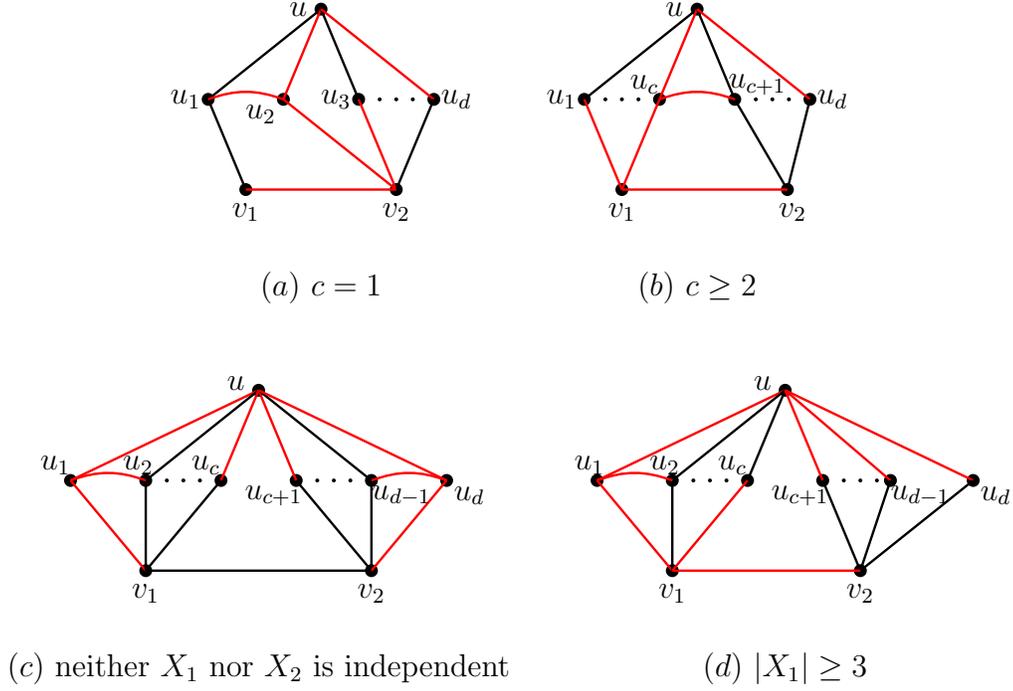
			\node at(8,-6.5) {Figure 2. HISTs (in red) in the proof of Fact 3 };
			
		\end{tikzpicture}
	\end{center}	
	%---------------
	
	\vskip -2mm
	
	\noindent{\bf Fact 3.} We can assume that $d< n-3$.\medskip
	
	\noindent{\bf Proof of Fact 3.}	Suppose that $d= n-3$. Then $N(u)=X_1\cup X_2$ and $d\ge 6$ as $n\ge 8$. If $E(X_1,X_2)\neq \emptyset$, say $u_cu_{c+1}\in E(G)$, then $G$ has a HIST $T$ (see Figure 2 (a) and (b)) with 
	\begin{eqnarray}
		E(T)=\left\{
		\begin{array}{ll}
			E(u,N(u)\backslash \{u_1,u_3\})\bigcup \{u_1u_2,u_2v_2,v_1v_2,v_2u_3\},&\mbox{if}~ c=1,\\ \nonumber			E(u,N(u)\backslash \{u_1,u_{c+1}\})\bigcup \{u_1v_1,u_cv_1,u_cu_{c+1},v_1v_2\},  &\mbox{if}~ c\ge 2.\nonumber
		\end{array}
		\right.
	\end{eqnarray}
	So we can assume that $E(X_1,X_2)= \emptyset$. First we note that $N(u)=X_1\cup X_2$ is not an independent set, for otherwise, one has $$2|E(G)|=2(1+2(n-3))=4n-10<n^2-7n+18$$ as $n\ge 8$, a contradiction with \eqref{q4}. If neither $X_1$ nor $X_2$ is independent, say $u_1u_2\in E(G)$ and $u_{d-1}u_d\in E(G)$, then $G$ has a HIST $T$ with 
	
	~~~~$E(T)=E(u,N(u)\setminus\{u_2,u_{d-1}\})\bigcup\{u_1u_{2},u_1v_1,u_{d-1}u_d,u_dv_2\}$ (see Figure 2 (c)). 
	
	\noindent Thus, we can assume that $X_1$ is not independent and $X_2$ is independent. Then $|X_1|\ge 2$. If $|X_1|= 2$, then $2|E(G)|=4n-8<n^2-7n+18$ as $n\ge 8$, a contradiction with \eqref{q4}. Hence, $|X_1|\ge 3$. Without loss of generality, we can assume $u_1u_2\in E(G)$, then $G$ has a HIST $T$ with $E(T)=E(u,N(u)\setminus\{u_2,u_{c}\})\bigcup\{u_1u_{2},u_1v_1,v_1u_c,v_1v_2\}$ (see Figure 2 (d)).\q
	
	\begin{center}
		\begin{tikzpicture}[line width =0.9pt]	
			\filldraw (12,-1.8) circle (.07);
			
			\filldraw (9,-3) circle (.07);
			\filldraw (9.28,-3) circle (.01);		
			\filldraw (9.5,-3) circle (.01);
			\filldraw (9.72,-3) circle (.01);
			\filldraw (10,-3) circle (.07);
			\filldraw (11,-3) circle (.07);
			\filldraw (11.28,-3) circle (.01);		
			\filldraw (11.5,-3) circle (.01);
			\filldraw (11.72,-3) circle (.01);
			\filldraw (12,-3) circle (.07);
			\filldraw (13,-3) circle (.07);
			\filldraw (14,-3) circle (.07);
			\filldraw (14.28,-3) circle (.01);		
			\filldraw (14.5,-3) circle (.01);
			\filldraw (14.72,-3) circle (.01);
			\filldraw (15,-3) circle (.07);
			
			\filldraw (9.5,-4.2) circle (.07);
			\filldraw (12,-4.2) circle (.07);
			
			\draw [red](12,-1.8) -- (9,-3); 
			\draw [red](12,-1.8) -- (10,-3); 
			\draw [red](12,-1.8) -- (11,-3); 
			\draw (12,-1.8) -- (12,-3); 
			\draw [red](12,-1.8) -- (13,-3); 
			\draw (12,-1.8) -- (14,-3); 
			\draw [red](12,-1.8) -- (15,-3); 
			
			\draw (9.5,-4.2) -- (9,-3); 
			\draw (9.5,-4.2) -- (10,-3);
			
			\draw (12,-4.2) -- (11,-3); 
			\draw [red](12,-4.2) -- (12,-3); 
			\draw [red](12,-4.2) -- (13,-3); 
			
			\draw [red](12,-4.2) -- (9.5,-4.2); 
			
			\draw[red] (13,-3) to [bend left=20] (14,-3);
			
			\node at(11.8,-1.6) {$u$};
			\node at(8.8,-3.2) {$u_1$};
			\node at(9.7,-3.2) {$u_c$};
			\node at(10.8,-3.3) {$u_{c+1}$};
			\node at(12.4,-2.9) {$u_{d-1}$};
			\node at(13.2,-3.3) {$u_d$};
			\node at(14.2,-3.3) {$u_j$};
			\node at(15.6,-3.3) {$u_{\Delta(G)}$};
			
			\node at(9.5,-4.5) {$v_1$};
			\node at(12,-4.5) {$v_2$};
			
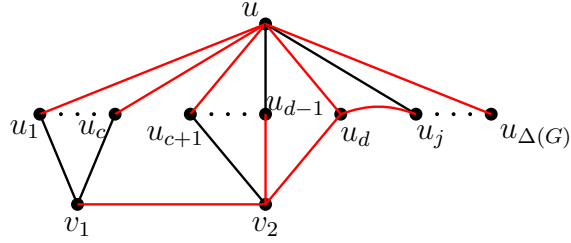
			\node at(12.5,-5.5) {Figure 3. HISTs (in red) when $X\not=\emptyset$ , $d-c\ge2$ };	
			
		\end{tikzpicture}
	\end{center}	
	
	By Fact 3, $X\not=\emptyset$. Since $G$ is 2-connected, $G-u$ is connected, and then $E(N(v_1)\cup N(v_2), X)\neq \emptyset$, which implies that exist some $u_i~(1\le i\le d)$ and $u_j~({d+1}\le j\le {\Delta(G)})$ such that  $u_iu_j\in E(G)$. Without loss of generality, we can assume that $i=d$. If $d-c\ge2$, then $G$ has a HIST $T$ with $E(T)=E(u,N(u)\backslash\{u_{d-1},u_j\})\bigcup\{u_du_j, u_{d-1}v_2,u_dv_2,v_1v_2\}$ (see Figure 3). 
	So we may assume $d-c=1$.
	
	\vskip.2cm
	\begin{center}
		\begin{tikzpicture}[line width =0.9pt]	
			\filldraw (3.5,-1.8) circle (.07);
			
			\filldraw (1,-3) circle (.07);
			\filldraw (1.28,-3) circle (.01);		
			\filldraw (1.5,-3) circle (.01);
			\filldraw (1.72,-3) circle (.01);
			
			\filldraw (2.28,-3) circle (.01);		
			\filldraw (2.5,-3) circle (.01);
			\filldraw (2.72,-3) circle (.01);
			
			\filldraw (2,-3) circle (.07);
			\filldraw (2.28,-3) circle (.01);		
			\filldraw (2.5,-3) circle (.01);
			\filldraw (2.72,-3) circle (.01);
			
			\filldraw (7,-3) circle (.07); 
			\filldraw (3,-3) circle (.07);
			\filldraw (4,-3) circle (.07);
			
			\filldraw (4.28,-3) circle (.01);		
			\filldraw (4.5,-3) circle (.01);
			\filldraw (4.72,-3) circle (.01);
			\filldraw (5,-3) circle (.07);
			
			\filldraw (5.28,-3) circle (.01);		
			\filldraw (5.5,-3) circle (.01);
			\filldraw (5.72,-3) circle (.01);
			\filldraw (7,-3) circle (.07); 
			\filldraw (6,-3) circle (.07);
			
			\filldraw (6.28,-3) circle (.01);		
			\filldraw (6.5,-3) circle (.01);
			\filldraw (6.72,-3) circle (.01);
			\filldraw (7,-3) circle (.07);
			
			\filldraw (2.5,-4.5) circle (.07);
			\filldraw (4,-4.5) circle (.07);
			
			\draw[red] (3.5,-1.8) -- (1,-3);
			\draw (3.5,-1.8) -- (3,-3);
			\draw [red](3.5,-1.8) -- (4,-3);
			\draw (3.5,-1.8) -- (5,-3);
			\draw [red](3.5,-1.8) -- (6,-3);
			\draw [red](3.5,-1.8) -- (7,-3);
			\draw [red](3.5,-1.8) -- (2,-3);
			
			\draw[red] (2.5,-4.5) -- (1,-3);
			\draw[red] (2.5,-4.5) -- (3,-3);
			\draw (4,-4.5) -- (4,-3);
			\draw[red] (2.5,-4.5) -- (4,-4.5);
			\draw (2.5,-4.5) -- (2,-3);
			
			\draw[red] (1,-3) to [bend left=20] (5,-3); 
			\draw (4,-3) to [bend right=20] (6,-3); 
			
			\node at(3.2,-1.8) {$u$};
			\node at(0.8,-3.2) {$u_1$};
			\node at(2.7,-3.2) {$u_c$};
			\node at(3.7,-3.2) {$u_d$};
			\node at(5,-2.8) {$u_p$};
			\node at(5.9,-3.3) {$u_j$};
			\node at(7.2,-3.3) {$u_{\Delta(G)}$};
			\node at(2.5,-4.8) {$v_1$};
			\node at(4,-4.8) {$v_2$};
			%%%%%%%%%%%%%%%%%%%%%%%%%%%%%%%%%%%%%%%%%%%%%%%%%%%%%%%%%%%5
			\filldraw (11,-1.8) circle (.07);
			
			\filldraw (9,-3) circle (.07);
			
			\filldraw (9.28,-3) circle (.01);		
			\filldraw (9.5,-3) circle (.01);
			\filldraw (9.72,-3) circle (.01);
			\filldraw (9,-3) circle (.07); 
			\filldraw (10,-3) circle (.07);
			
			\filldraw (10.28,-3) circle (.01);		
			\filldraw (10.5,-3) circle (.01);
			\filldraw (10.72,-3) circle (.01); 
			\filldraw (11,-3) circle (.07);
			\filldraw (12,-3) circle (.07);
			
			\filldraw (12.28,-3) circle (.01);		
			\filldraw (12.5,-3) circle (.01);
			\filldraw (12.72,-3) circle (.01);
			\filldraw (13,-3) circle (.07);
			
			\filldraw (14,-3) circle (.07);
			\filldraw (13.28,-3) circle (.01);		
			\filldraw (13.5,-3) circle (.01);
			\filldraw (13.72,-3) circle (.01);
			\filldraw (14.28,-3) circle (.01);		
			\filldraw (14.5,-3) circle (.01);
			\filldraw (14.72,-3) circle (.01);
			\filldraw (15,-3) circle (.07);
			
			\filldraw (10.5,-4.5) circle (.07);
			\filldraw (12,-4.5) circle (.07);
			
			\draw [red](11,-1.8) -- (9,-3);
			\draw (11,-1.8) -- (10,-3);
			\draw [red](11,-1.8) -- (11,-3);
			\draw [red] (11,-1.8) -- (12,-3);
			\draw [red] (11,-1.8) -- (13,-3);
			\draw  (11,-1.8) -- (14,-3);
			\draw [red] (11,-1.8) -- (15,-3);

			\draw[red](10.5,-4.5) -- (9,-3);
			\draw (10.5,-4.5) -- (10,-3);
			\draw (10.5,-4.5) -- (11,-3);
			
			\draw[red] (12,-4.5) -- (12,-3);
			
			\draw (12,-4.5) -- (10.5,-4.5);
			
			\draw[red] (9,-3) to [bend left=20] (10,-3); 
			\draw[red] (12,-3) to [bend left=20] (14,-3); 
			
			\node at(10.7,-1.8) {$u$};
			\node at(8.8,-3.3) {$u_1$};
			\node at(9.8,-3.3) {$u_p$};
			\node at(10.65,-3.3) {$u_c$};
			\node at(11.7,-3.3) {$u_d$};
			\node at(13.9,-3.3) {$u_j$};
			\node at(15.2,-3.3) {$u_{\Delta(G)}$};
			\node at(10.5,-4.8) {$v_1$};
			\node at(12,-4.8) {$v_2$}; 
			
			\node at(4.5,-5.5) {$(a)~u_p\in X_2\cup X$};
			\node at(11.2,-5.5) {$(b)~u_p\in X_1$};
			
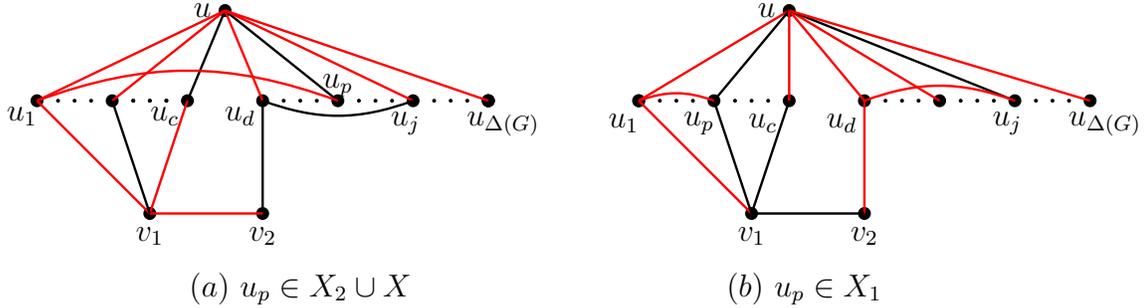
			\node at(8,-6.5) {Figure 4. HISTs (in red) when $d-c=1$ };
		\end{tikzpicture}
	\end{center}	
	
	On the other hand, we can assume that $d(u_i)=2$ for all $1\le i \le c$. Otherwise, we can assume that $d(u_1)>2$. Then $N(u_1)\setminus \{u,v_1,v_2\}\neq \emptyset$. Let $u_p\in N(u_1)\setminus \{u,v_1,v_2\}$. Then $G$ has a HIST $T$ (see Figure 4 (a) and (b)) with 
	\begin{eqnarray}
		E(T)=\left\{
		\begin{array}{ll}
			E(u,N(u)\backslash\{u_{c},u_p\})\bigcup\{u_1u_p, u_{1}v_1, v_1v_2, v_1u_c\},&\mbox{if}~ u_p\in X_2\cup X,\\ \nonumber
			E(u,N(u)\backslash\{u_{p},u_j\})\bigcup\{u_1u_p, u_{1}v_1, u_dv_2,u_du_j\},  &\mbox{if}~ u_p\in X_1.\nonumber
		\end{array}
		\right.
	\end{eqnarray}

	\noindent{\bf Fact 4.} $c=1$.\medskip
	
	\noindent{\bf Proof of Fact 4.} Note that $d(v_1)=c+1$, $d(v_2)=2$, $d(u_{c+1})\le n-c-2$ and $ d(u_i)\le n-3-c $  for $c+2 \le i  \le n-3$. Suppose that $c\ge 2$. Then 		
	\begin{eqnarray*}\label{e1}
		2|E(G)|&=&d(u)+d(v_1)+d(v_2)+d(u_{c+1})+\sum_{i=1}^{c}d(u_i)+\sum_{j=c+2}^{n-3}d(u_j)\nonumber\\
		&\leq& (n-3)+(c+1)+2+(n-2-c)+2c+(n-4-c)(n-3-c)\nonumber\\ 
		&=& c^2-(2n-9)c+n^2-5n+10\nonumber\\
		&\le&
		\max\{2^2-2(2n-9)+n^2-5n+10,(n-4)^2-(n-4)(2n-9)n^2-5n+10\}\nonumber\\ 
		&=& n^2-9n+32
		< n^2-7n+18,
	\end{eqnarray*}
	where the last inequality follows from $n\ge 8$, a contradiction with \eqref{q4}.  \q \medskip

	By Fact 4, we have $c=1$, then $d(u_1)=d(v_1)=d(v_2)=2$. Furthermore, $G$ is a subgraph of $B_n$, and then by Lemma \ref{le1} and the assumption, we obtain $G\cong B_n$. 
	
	\vskip.6cm
	
	\begin{center}
		\begin{tikzpicture}[line width =0.9pt]			
			
			\filldraw (1.5,-1.8) circle (.07);	
			\filldraw (0,-3) circle (.07);
			\filldraw (-1,-3) circle (.07);	
			\filldraw (4,-3) circle (.07);	
			\filldraw (1,-3) circle (.07);	
			
			\filldraw (2,-3) circle (.07);
			\filldraw (2.28,-3) circle (.01);		
			\filldraw (2.5,-3) circle (.01);
			\filldraw (2.72,-3) circle (.01);
			\filldraw (3,-3) circle (.07);	
			
			\filldraw (0.5,-4.2) circle (.07);	
			\filldraw (2.5,-4.2) circle (.07);	
			
			\draw (1.5,-1.8) -- (0,-3);
			\draw[red](1.5,-1.8) -- (-1,-3);
			\draw[red](1.5,-1.8) -- (4,-3);
			
			\draw [red](1.5,-1.8) -- (1,-3);
			\draw [red](1.5,-1.8) -- (2,-3);
			\draw (1.5,-1.8) -- (3,-3);
			
			\draw (0.5,-4.2) -- (-1,-3);
			\draw (2.5,-4.2) -- (4,-3);
			
			\draw(0.5,-4.2) -- (0,-3);
			\draw[red] (0.5,-4.2) -- (1,-3);
			
			\draw [red](2.5,-4.2) -- (2,-3);
			\draw (2.5,-4.2) -- (3,-3);
			
			%		\draw [red](0.5,-4.2) -- (2.5,-4.2);
			\draw [red](0,-3) to [bend left=20] (2,-3);
			\draw [red](1,-3) to [bend right=20] (3,-3); 
			
			\node at(1.2,-1.7) {$u$};
			\node at(-0.3,-3) {$u_i$};
			\node at(0.7,-3.15) {$u_j$};
			\node at(1.7,-3) {$u_k$};
			\node at(3.3,-3) {$u_l$};
			\node at(0.5,-4.5) {$v_1$};
			\node at(2.5,-4.5) {$v_2$};
			
			\node at(1.5,-5.5) {$(a)~i\not=j,k\not=l$};
			
			%%%%%%%%%%%%%%%%%%%%%%%%%%%%%%%%%%%%%%%%%%%	
			
			\filldraw (6.5,-1.8) circle (.07);
			
			\filldraw (5,-3) circle (.07);
			\filldraw (5.28,-3) circle (.01);		
			\filldraw (5.5,-3) circle (.01);
			\filldraw (5.72,-3) circle (.01);
			\filldraw (6,-3) circle (.07);
			
			\filldraw (7,-3) circle (.07);
			\filldraw (7.28,-3) circle (.01);		
			\filldraw (7.5,-3) circle (.01);
			\filldraw (7.72,-3) circle (.01);
			\filldraw (8,-3) circle (.07);
			\filldraw (9,-3) circle (.07);
			
			\filldraw (5.5,-4.2) circle (.07);
			\filldraw (8,-4.2) circle (.07);
			
			\draw (6.5,-1.8) -- (6,-3);
			\draw[red](6.5,-1.8) -- (5,-3);
			\draw[red](6.5,-1.8) -- (7,-3);
			\draw[red](6.5,-1.8) -- (9,-3);
			\draw(6.5,-1.8) -- (8,-3);
			
			\draw(5.5,-4.2) -- (5,-3);
			\draw [red](5.5,-4.2) -- (6,-3);
			\draw[red](8,-4.2) -- (7,-3);
			\draw (8,-4.2) -- (8,-3);
			\draw (8,-4.2) -- (9,-3);

			%		\draw [red](5.5,-4.2) -- (7.7,-4.2);
			
			\draw[red] (6,-3) to [bend left=20] (7,-3); 
			\draw[red] (6,-3) to [bend right=20] (8,-3);
			
			\node at(6.2,-1.7) {$u$};
			%		\node at(5,-2.8) {$u_1$};
			\node at(5.8,-2.8) {$u_j$};
			\node at(7.3,-2.8) {$u_k$};
			\node at(8.3,-3) {$u_l$};
			\node at(5.5,-4.5) {$v_1$};
			\node at(8,-4.5) {$v_2$};
			
			\node at(6.5,-5.5) {$(b)~i=j,k\not=l$};
			
			\node at(6.5,-6.8) {Figure 5. HISTs are in red};	
			
			%%%%%%%%%%%%%%%%%%%%%%%%%%%%%%%%%%%%%%%%%%%	
			
			\filldraw (12.5,-1.8) circle (.07);
			
			\filldraw (10,-3) circle (.07);
			\filldraw (11,-3) circle (.07);
			\filldraw (11.28,-3) circle (.01);		
			\filldraw (11.5,-3) circle (.01);
			\filldraw (11.72,-3) circle (.01);
			\filldraw (12,-3) circle (.07);
			
			\filldraw (13,-3) circle (.07);
			\filldraw (13.28,-3) circle (.01);		
			\filldraw (13.5,-3) circle (.01);
			\filldraw (13.72,-3) circle (.01);
			\filldraw (14,-3) circle (.07);
			%		\filldraw (15,-3) circle (.07);
			
			\filldraw (11,-4.2) circle (.07); 
			\filldraw (13.5,-4.2) circle (.07); 
			
			\draw [red](12.5,-1.8) -- (10,-3);
			\draw[red](12.5,-1.8) -- (11,-3);
			\draw(12.5,-1.8) -- (12,-3);
			\draw(12.5,-1.8) -- (13,-3);
			\draw (12.5,-1.8) -- (14,-3);
			%		\draw [red](12.5,-1.8) -- (15,-3);
			
			\draw(11,-4.2) -- (10,-3);
			\draw[red](11,-4.2) -- (11,-3);
			\draw (11,-4.2) -- (12,-3);
			
			\draw[red](13.5,-4.2) -- (13,-3);
			\draw (13.5,-4.2) -- (14,-3);
			%		\draw [red](14,-4.2) -- (15,-3);
			
			%		\draw (11,-4.2) -- (14,-4.2);
			
			\draw[red] (11,-3) to [bend right=20] (13,-3); 
			\draw[red] (12,-3) to [bend left=20] (13,-3); 
			
			\node at(12.2,-1.7) {$u$};
			%		\node at(9.8,-2.8) {$u_1$};
			\node at(10.8,-2.8) {$u_i$};
			\node at(11.8,-2.8) {$u_j$};
			\node at(13.2,-2.8) {$u_l$};
			%		\node at(13.8,-2.8) {$u_{d-1}$};
			%		\node at(15.3,-3) {$u_d$};
			\node at(11,-4.5) {$v_1$};
			\node at(13.5,-4.5) {$v_2$};
			
			\node at(12,-5.5) {$(c)~i\not=j,k=l$};
			%%%%%%%%%%%%%%%%%%%%%%%%%%%%%%%%%%%%%%%%%%%	
			
		\end{tikzpicture}
	\end{center}
	
	\vskip.2cm

	\noindent {\bf Subcase 2.2.} $v_1v_2\notin E(G)$.
	
	\vskip.2cm
	
	In this subcase, we denote $N(v_1)=\{u_1,\ldots,u_p\}$ and $N(v_2)=\{u_{p+1},\ldots,u_q\}$. Since $G$ is $2$-connected, we obtain $p\ge2$, $q-p\ge2$ and $q\le n-3$. Denote $Y=N(u)\backslash (N(v_1)\cup N(v_2))$. Since $G-u$ is connected, we have 
	\begin{equation} \label{Th2.2-e1} 
		E(N(v_1),N(v_2)\cup Y)\neq \emptyset,
	\end{equation} 
	\begin{equation} \label{Th2.2-e2} 
		E(N(v_2),N(v_1)\cup Y)\neq \emptyset.
	\end{equation} 
	
	\vskip.4cm

	\noindent{\bf Fact 5.} We can assume that $q< n -3$.\medskip
	
	\noindent{\bf Proof of Fact 5.}	Suppose that $q=n-3$, then by \eqref{Th2.2-e1}, one obtains $E(N(v_1),N(v_2))\neq \emptyset$. 
	If $|E(N(v_1),N(v_2))|\ge2$, then there exist some $u_i,u_j\in N(v_1)$ and $u_k,u_l\in N(v_2)$ such that $u_iu_k,u_ju_l\in E(G)$, and thus $G$ has a HIST $T$ (see Figure 5 (a)-(c)) with 
	\begin{eqnarray}
		E(T)=\left\{
		\begin{array}{ll}
			E(u,N(u)\backslash \{u_i,u_l\})\bigcup \{u_jv_1,u_ju_l,u_iu_k,u_kv_2\}, &\mbox{if}~ u_i\not=u_j, u_k\not=u_l,\\ \nonumber
			E(u,N(u) \backslash \{u_j,u_l\})\bigcup \{u_jv_1,u_ju_l,u_ju_k,u_kv_2\}, &\mbox{if}~ u_i=u_j, u_k\not=u_l,\\ \nonumber
			E(u,N(u)\backslash \{u_j,u_l\})\bigcup \{u_iv_1,u_iu_l,u_ju_l,u_lv_2\},  &\mbox{if}~ u_i\not=u_j, u_k=u_l.\nonumber
		\end{array}
		\right.
	\end{eqnarray}
	So $|E(N(v_1),N(v_2))|=1$. Without loss of generality, we assume $u_p, u_{p+1}\in E(G)$.  If both $N(v_1)$ and $N(v_2)$ are independent, then $2|E(G)|=2(1+2(n-3))=4n-10<n^2-7n+18$ as $n\ge 8$, a contradiction with \eqref{q4}. Therefore either $N(v_1)$ or $N(v_2)$ is not independent. Without loss of generality, we can let $u_s, u_t\in N(v_1) $ with $u_su_t\in E(G)$. Then $G$ has a HIST $T$ with $E(T)=E(u,N(u))\backslash \{u_s,u_p\})\bigcup \{u_tv_1,u_su_t,u_pu_{p+1},u_{p+1}v_2\}$ (see Figure 6 (a)). \q
	
	\vskip.2cm
	
	\begin{center}
		\begin{tikzpicture}[line width =0.9pt]			
			
			\filldraw (1.5,-1.8) circle (.07);	
			\filldraw (0,-3) circle (.07);
			\filldraw (-1,-3) circle (.07);	
			%		\filldraw (4,-3) circle (.07);	
			\filldraw (1,-3) circle (.07);	
			
			\filldraw (2,-3) circle (.07);
			\filldraw (2.28,-3) circle (.01);		
			\filldraw (2.5,-3) circle (.01);
			\filldraw (2.72,-3) circle (.01);
			\filldraw (3,-3) circle (.07);	
			
			\filldraw (0.5,-4.2) circle (.07);	
			\filldraw (2.5,-4.2) circle (.07);	
			
			\draw[red](1.5,-1.8) -- (0,-3);
			\draw(1.5,-1.8) -- (-1,-3);
			\draw[red](0,-3) -- (-1,-3);
			
			\draw(1.5,-1.8) -- (1,-3);
			\draw [red](1.5,-1.8) -- (2,-3);
			\draw[red](1.5,-1.8) -- (3,-3);
			
			\draw (0.5,-4.2) -- (-1,-3);
			%		\draw (2.5,-4.2) -- (4,-3);
			
			\draw[red](0.5,-4.2) -- (0,-3);
			\draw(0.5,-4.2) -- (1,-3);
			
			\draw [red](2.5,-4.2) -- (2,-3);
			\draw (2.5,-4.2) -- (3,-3);
			
			%		\draw [red](0.5,-4.2) -- (2.5,-4.2);
			\draw [red](1,-3) to [bend left=20] (2,-3);
			%		\draw [red](1,-3) to [bend right=20] (3,-3); 
			
			\node at(1.2,-1.7) {$u$};
			\node at(-1.3,-3) {$u_s$};
			\node at(-0.3,-3.2) {$u_t$};
			\node at(0.7,-2.8) {$u_p$};
			\node at(1.7,-3.3) {$u_{p+1}$};
			\node at(3.3,-3) {$u_q$};
			\node at(0.5,-4.5) {$v_1$};
			\node at(2.5,-4.5) {$v_2$};
			
			\node at(1.5,-5.5) {$(a)$};
			
			%%%%%%%%%%%%%%%%%%%%%%%%%%%%%%%%%%%%%%%%%%%	
			
			\filldraw (5.5,-1.8) circle (.07);
			
			\filldraw (4,-3) circle (.07);
			\filldraw (4.28,-3) circle (.01);		
			\filldraw (4.5,-3) circle (.01);
			\filldraw (4.72,-3) circle (.01);
			\filldraw (5,-3) circle (.07);
			
			\filldraw (6,-3) circle (.07);
			\filldraw (6.28,-3) circle (.01);		
			\filldraw (6.5,-3) circle (.01);
			\filldraw (6.72,-3) circle (.01);
			\filldraw (7,-3) circle (.07);
			\filldraw (8,-3) circle (.07);
			
			\filldraw (4.5,-4.2) circle (.07);
			\filldraw (7,-4.2) circle (.07);
			
			\draw (5.5,-1.8) -- (5,-3);
			\draw[red](5.5,-1.8) -- (4,-3);
			\draw[red](5.5,-1.8) -- (6,-3);
			\draw(5.5,-1.8) -- (8,-3);
			\draw[red](5.5,-1.8) -- (7,-3);
			
			\draw[red](4.5,-4.2) -- (4,-3);
			\draw(4.5,-4.2) -- (5,-3);
			\draw[red](7,-4.2) -- (6,-3);
			\draw (7,-4.2) -- (7,-3);
			%		\draw (8,-4.2) -- (9,-3);

			%		\draw [red](5.5,-4.2) -- (7.7,-4.2);
			
			\draw[red] (5,-3) to [bend left=20] (6,-3); 
			\draw[red] (4,-3) to [bend right=20] (8,-3);
			
			\node at(5.2,-1.7) {$u$};
			\node at(4,-2.8) {$u_i$};
			\node at(4.8,-2.8) {$u_p$};
			\node at(5.75,-3.2) {$u_{p+1}$};
			\node at(7.3,-3) {$u_q$};
			\node at(8.3,-3) {$u_j$};
			\node at(4.5,-4.5) {$v_1$};
			\node at(7,-4.5) {$v_2$};
			
			\node at(6,-5.5) {$(b)$};
			
			%\node at(6.5,-6.5) {Figure 5. (a) HIST when $Y=\emptyset$, (b) HIST when $Y\not=\emptyset$ and $E(N(v_1),N(v_2))\neq \emptyset$,}; 
			%\node at(6.5,-7) { (c) HIST when $Y\not=\emptyset$ and $E(N(v_1),N(v_2))=\emptyset$ (in red) };	
			\node at(6.5,-6.5) {Figure 6. HISTs are in red };	
			
			%%%%%%%%%%%%%%%%%%%%%%%%%%%%%%%%%%%%%%%%%%%	
			
			\filldraw (11.5,-1.8) circle (.07);
			
			\filldraw (9,-3) circle (.07);
			\filldraw (10,-3) circle (.07);
			\filldraw (9.28,-3) circle (.01);		
			\filldraw (9.5,-3) circle (.01);
			\filldraw (9.72,-3) circle (.01);
			\filldraw (11,-3) circle (.07);
			
			\filldraw (12,-3) circle (.07);
			\filldraw (11.28,-3) circle (.01);		
			\filldraw (11.5,-3) circle (.01);
			\filldraw (11.72,-3) circle (.01);
			\filldraw (13,-3) circle (.07);
			\filldraw (14,-3) circle (.07);
			
			\filldraw (9.5,-4.2) circle (.07); 
			\filldraw (11.5,-4.2) circle (.07); 
			
			\draw [red](11.5,-1.8) -- (9,-3);
			\draw[red](11.5,-1.8) -- (10,-3);
			\draw[red](11.5,-1.8) -- (11,-3);
			\draw[red](11.5,-1.8) -- (12,-3);
			\draw (11.5,-1.8) -- (13,-3);
			\draw(11.5,-1.8) -- (14,-3);
			
			\draw[red](9.5,-4.2) -- (9,-3);
			\draw(9.5,-4.2) -- (10,-3);
			%		\draw (10,-4.2) -- (11,-3);
			
			\draw[red](11.5,-4.2) -- (12,-3);
			\draw(11.5,-4.2) -- (11,-3);
			%		\draw (12.5,-4.2) -- (13,-3);
			%		\draw [red](14,-4.2) -- (15,-3);
			
			%		\draw (11,-4.2) -- (14,-4.2);
			
			\draw[red] (14,-3) to [bend right=20] (12,-3); 
			\draw[red] (9,-3) to [bend right=20] (13,-3); 
			
			\node at(11.2,-1.7) {$u$};
			\node at(8.8,-2.8) {$u_1$};
			\node at(9.8,-2.8) {$u_p$};
			\node at(10.7,-3.2) {$u_{p+1}$};
			\node at(12.2,-2.75) {$u_q$};
			\node at(13,-3.3) {$u_j$};
			\node at(14,-3.3) {$u_k$};
			\node at(9.5,-4.5) {$v_1$};
			\node at(11.5,-4.5) {$v_2$};
			
			\node at(11.5,-5.5) {$(c)$};
			%%%%%%%%%%%%%%%%%%%%%%%%%%%%%%%%%%%%%%%%%%%	

		\end{tikzpicture}
	\end{center}
	
	By Fact 5, $Y\not= \emptyset$. If $E(N(v_1),N(v_2))\neq \emptyset$, say $u_pu_{p+1}\in E(G)$, then $E(N(v_1)\cup N(v_2),Y)\neq \emptyset$ as $G-u$ is connected. Assume that there exists some $u_i\in N(v_1)$ and $u_j\in Y $ such that $u_iu_j\in E(G)$, then $G$ contains a HIST $T$ with $E(T)=E(u,N(u))\backslash \{u_p,u_j\})\medcup\{u_iv_1,u_iu_j,u_pu_{p+1},u_{p+1}v_2\}$ (see Figure 6 (b)).\medskip
	
	If $E(N(v_1),N(v_2))=\emptyset$, then by \eqref{Th2.2-e1} and \eqref{Th2.2-e2},  $E(N(v_a),Y)\not=\emptyset$ for $a=1,2$. Assume that $u_1u_j\in E(G)$ and  $u_qu_k\in E(G)$ for some $q+1\leq j,k\leq n-3$. If $j\neq k$, then $G$ contains a HIST $T$ with $E(T)=E(u,N(u)\backslash \{u_j,u_k\})\bigcup \{u_1v_1,u_1u_j,u_qv_2,u_qu_k\}$ (see Figure 6 (c)). So we can assume that $(N(u_l)\backslash \{u,v_1,v_2\})\cap Y=\{u_j\}$ for $1\le l\leq q$, where $u_j\in Y$. It follows that $d(v_1)=p, ~d(v_2)=q-p, ~d(u_j)\leq n-3, ~d(u_g)\le 3 $ for $1 \le g \le q$, $d(u_h)\le n-q-3 $ for $h\neq j,q+1 \le h\le n-3$. Then 		
	\begin{eqnarray*}\label{e2}
		2|E(G)|&=&d(u)+d(v_1)+d(v_2)+\sum_{g=1}^{q}d(u_g)+\sum_{h=q+1}^{n-3}d(u_h)\nonumber\\ 
		&\leq& (n-3)+p+(q-p)+(n-3)+3q+(n-q-3)(n-q-4) \nonumber\\ 
		%&=& 4q+2(n-3)+(n-q-3)(n-q-4)\nonumber\\
		&=& q^2-(2n-11)q+n^2-5n+6\nonumber\\   	
		&\leq& 
		\max \{4^2-4(2n-11)+n^2-5n+6,\nonumber\\ 
		&&~~~~~~~~(n-4)^2-(2n-11)(n-4)+n^2-5n+6\}\nonumber\\ 
		&=&  n^2-13n+66 
		\le  n^2-7n+18,
	\end{eqnarray*}
	where the last inequality follows from $n\ge 8$, a contradiction with \eqref{q4}.

	Therefore, the proof of Theorem \ref{Thm2} is complete.	\q
	\vskip-2mm
	
	\section*{Declaration of Competing Interest}
	
	The authors declare that they have no known competing financial interests or personal relationships that could have appeared to influence the work reported in this paper.


\begin{thebibliography}{99}
		
		\bibitem{Michael1990} M.O. Albertson, D.M. Berman, J.P. Hutchinson, C. Thomassen, Graphs with homeomorphically irreducible spanning trees, J. Graph Theory 14 (1990) 247-258.
		
		\bibitem{Babai2009}	L. Babai, B. Guiduli, Spectral extrema for graphs: the Zarankiewicz problem, Electron. J. Comb. 16 (2009) \#R123.
		
		\bibitem{Chen2019} M.Z. Chen, A.M. Liu, X.D. Zhang, Spectral extremal results with forbidding linear forests,
		Graphs Comb. 35 (2019) 335-351.
		
		\bibitem{CRS10} D. Cvetkovi\'c, P. Rowlinson, S. Simi\'c, An Introduction to the Theory of Graph Spectra, Cambridge University Press, Cambridge, 2010.
		
		\bibitem{Fiedler2010} M. Fiedler, V. Nikiforov, Spectral radius and Hamiltonicity of graphs, Linear Algebra Appl. 432 (9) (2010) 2170-2173.
		
		\bibitem{Furuya2025} M. Furuya, A. Saito, S. Tsuchiya, Refinements of degree conditions for the existence of a spanning tree without small degree stems, Discrete Math. 348 (2) (2025) 114307.
		
		\bibitem{Hill1974} A. Hill, Graphs with homeomorphically irreducible spanning trees, LMS Lect. Note Ser., 13 (1974) 61-68.
		
		\bibitem{Hong1988} Y. Hong, A bound on the spectral radius of graphs, Linear Algebra Appl. 108 (1988) 135-139.
		
		\bibitem{Ito2022} T. Ito, S. Tsuchiya, Degree sum conditions for the existence of homeomorphically irreducible spanning trees, J. Graph Theory 99 (1) (2022) 162-170.
		
		\bibitem{Hong1999} Y. Hong, J.L Shu, A Sharp Upper Bound of the Spectral Radius of Graphs, J. Comb. Theory, Ser. B 81(2001)  177-183.
		
		\bibitem{Li2024} Y. Li, F.M. Dong, X.L. Hu, H.Q. Liu, A neighborhood union condition for the existence of a spanning tree without degree $2$ vertices, arxiv: 2412.07128.	
		
		%		\bibitem{Lu2012} M. Lu, H. Liu, F. Tian, Spectral radius and Hamiltonian graphs, Linear Algebra
		%		Appl. 437 (7) (2012) 1670-1674.
		
		\bibitem{Lu2024} J. Lu, L. Lu, Y. Li, Spectral radius of graphs forbidden $C_7$ or $C_6^{\Delta}$, Discrete Math. 347 (2)
		(2024) 113781.
		
		\bibitem{Min2022} G. Min, Z. Lou, Q. Huang, A sharp upper bound on the spectral radius of $C_5$-free/$C_6$-free
		graphs with given size, Linear Algebra Appl. 640 (2022) 162-78.
		
		%		\bibitem{Ning2015} B. Ning, J. Ge, Spectral radius and Hamiltonian properties of graphs, Linear Multilinear Algebra 63 (8) (2014) 1520-1530.
		
		\bibitem{Niki} V. Nikiforov, The spectral radius of graphs without paths and cycles of specified length, Linear Algebra Appl. 432 (2010) 2243-2256.
		
		\bibitem{Niki2007} V. Nikiforov, Bounds on graph eigenvalues II, Linear Algebra Appl. 427 (2007) 183-189.
		
		\bibitem{Wilf} H. Wilf, Spectral bounds for the clique and independence numbers of graphs, J. Comb. Theory, Ser. B 40 (1986) 113-117.
		
		\bibitem{Li2022} Y. Li, Y, Peng, The maximum spectral radius of non-bipartite graphs forbidding short odd
		cycles, Electron. J. Comb. 29 (4) (2022) \#P4.2.
		
		\bibitem{Zhai2012} M. Zhai, B. Wang, Proof of a conjecture on the spectral radius of $C_4$-free graphs, Linear Algebra
		Appl. 437 (2012) 1641-1647.
		
		\bibitem{Zhai2020} M. Zhai, H. Lin, Spectral extrema of graphs: forbidden hexagon, Discrete Math. 343 (10) (2020) 112028.
		
		\bibitem{Zhai2022J} M. Zhai, H. Lin, Spectral extrema of $K_{s,t}$-minor free graphs-On a conjecture of M. Tait., J. Comb. Theory, Ser. B 157 (2022) 184-215.
		
		\bibitem{Zhai2022E} M. Zhai, R. Liu, J. Xue, A unique characterization of spectral extrema for friendship graphs, Electron. J. Comb. 29 (3) (2022) \#P3.32.
		
		
		
		
		
		
		
		
		
		
	\end{thebibliography}
\end{document}